\documentclass[11pt]{amsart}
\usepackage{amsmath,amssymb,amsthm}
\usepackage{epsf}
\usepackage{epsfig}

\newcommand{\R}{\mathbb R}
\newcommand{\C}{\mathbb C}

\newtheorem{thm}{Theorem}[section]

\theoremstyle{definition}

\theoremstyle{remark}

\newcommand{\ds}{\displaystyle}

\begin{document}

\title[EXPLICIT SOLVING OF THE SYSTEM OF NATURAL PDE'S]
{Explicit Solving of the System of Natural PDE's\\ of Minimal Surfaces\\ in the Four-Dimensional Euclidean Space}%

\author{Georgi Ganchev and Krasimir Kanchev}

\address{Bulgarian Academy of Sciences, Institute of Mathematics and Informatics,
Acad. G. Bonchev Str. bl. 8, 1113 Sofia, Bulgaria}
\email{ganchev@math.bas.bg}%

\address {Department of Mathematics and Informatics, Todor Kableshkov University of Transport,
158 Geo Milev Str., 1574 Sofia, Bulgaria}%
\email{kbkanchev@yahoo.com}%


\begin{abstract}

The fact that minimal surfaces in the four-dimensional Euclidean space admit natural parameters
implies that any minimal surface is determined uniquely up to a motion by two
curvature functions, satisfying a system of two PDE's (the system of natural PDE's). In fact
this solves the problem of Lund-Regge for minimal surfaces. Using the corresponding result for
minimal surfaces in the three-dimensional Euclidean space, we solve explicitly the system of
natural PDE's, expressing any solution by virtue of two holomorphic functions in the Gauss plane.
We find the relation between two pairs of holomorphic functions (i.e. the class of pairs of holomorphic
functions) generating one and the same solution of the system of natural PDE's.

\end{abstract}

\maketitle

\thispagestyle{empty}

\section{Introduction}

Studying minimal surfaces in the four-dimensional Euclidean space $\R^4$,
Itoh proved in \cite{Itoh} that any minimal non-superconformal surface $M^2$ admits locally
special isothermal parameters. On the base of this result, de Azevero Tribuzy and
Guadalupe proved in \cite{T-G} the following theorem:

The Gauss curvature $K$ and the curvature of the normal connection $\varkappa$ (the normal curvature)
of a minimal non-superconformal surface, parameterized by special isothermal parameters, satisfy
the following system of  partial differential equations:
\begin{equation}\label{Eq-1}
\begin{array}{l}
\ds{(K^2 - \varkappa^2)^{\frac{1}{4}}\, \Delta \ln |\varkappa - K|} = 2(2K - \varkappa);\\
[2mm]
\ds{(K^2 - \varkappa^2)^{\frac{1}{4}}\, \Delta \ln |\varkappa + K| = 2(2K + \varkappa)}.
\end{array}
\end{equation}

Conversely, any solution ($K$, $\varkappa$) to system \eqref{Eq-1} determines uniquely (up to a
motion in $\R^4$) a minimal non-superconformal surface with Gauss
curvature $K$ and normal curvature $\varkappa$.

Further we call system \eqref{Eq-1} {\it the system of natural PDE's of minimal surfaces in $\R^4$}
and our aim is to solve explicitly this system.

All considerations in the paper are local.

Introducing natural parameters on any minimal surface in $\R^4$ reduces the number of the invariants
determining the surface to two: $K$ and $\varkappa$. Further, these two invariants satisfy the system
of two natural PDE's and determine the minimal surface uniquely up to a motion. It is clear that the
number of the invariants and the number of the PDE's can not be reduced further. Therefore this
solves the problem of Lund-Regge \cite{L-R} for minimal surfaces in $\R^4$.

The basic theorem in this paper is the following statement

{\bf Theorem 1. (Explicit solving of the system of natural PDE's of minimal surfaces)}
Let $K$ and $\varkappa$ be solutions of the system \eqref{Eq-2}. Then we have locally
\begin{equation}\label{Eq-3}
\begin{array}{l}
\ds{K=\frac{-8|w_1^\prime w_2^\prime|} {(|w_1|^2+1)(|w_2|^2+1)}\left(\frac{|w_1^\prime|^2}{(|w_1|^2+1)^2}+\frac{|w_2^\prime|^2}{(|w_2|^2+1)^2}\right),}\\
[6mm] \ds{\varkappa=\frac{8|w_1^\prime w_2^\prime|}
{(|w_1|^2+1)(|w_2|^2+1)}\left(\frac{|w_1^\prime|^2}{(|w_1|^2+1)^2}-\frac{|w_2^\prime|^2}{(|w_2|^2+1)^2}\right),}
\end{array}
\; w_1^\prime  w_2^\prime \neq 0,
\end{equation}
for some holomorphic functions $w_k$ \, $(k=1,2)$ \, in \, $\mathbb C$.

Conversely, any two functions $K$ and $\varkappa$, given by \eqref{Eq-3} ,
satisfy the system \eqref{Eq-2}.
\vskip 2mm
We shall say that the pair $(w_1, w_2)$ generates the solution \eqref{Eq-3}.
\vskip 2mm
There arises the following natural question: \emph{When two pairs of holomorphic functions
generate one and the same solution of \eqref{Eq-2}?}

The answer is given by the following statement

\textbf{Theorem 2.}
Let $(w_1, w_2)$ and $(\hat w_1, \hat w_2)$ be two pairs of holomorphic functions generating one and the same
solution $(K, \varkappa)$ of the system \eqref{Eq-2}.

Then
\begin{equation}\label{Eq-4}
\hat w_k = \frac{-\bar b_k +\bar a_k\,w_k}{a_k+b_k\,w_k},\quad a_k=\text{const}, \; b_k=\text{const},\; |a_k|^2+|b_k|^2=1; \; (k=1,2).
\end{equation}

Conversely, any two pairs of holomorphic functions related by \eqref{Eq-4} generate one and the same solution of \eqref{Eq-2}.
\vskip 2mm
The basic idea we use in solving the system of natural PDE's of minimal surfaces in $\R^4$ is to reduce it
to the solving of the natural equation of minimal surfaces in $\R^3$.

\section{Explicit solving of the system of natural PDE's}

In this section we find an explicit form of the solutions of the system \eqref{Eq-2}.

In \cite{G} the first author proved the following result:
\vskip 2mm
{Theorem A. (Explicit solving of the natural PDE of minimal surfaces in $\R^3$)
{\it Any solution $\nu > 0$ of the natural partial differential equation of minimal
surfaces

\begin{equation}\label{Eq-cite1}
\Delta \ln \nu + 2 \nu =0
\end{equation}
locally is given by the formula

\begin{equation}\label{Eq-cite2}
\nu =\frac{4|w^\prime|^2}{(|w|^2+1)^2}\,,\quad w^\prime \neq 0,
\end{equation}
where $w$ is a holomorphic function in $\C$.

Conversely, any function $\nu (x,y)$ of the type \eqref{Eq-cite2} is a solution to \eqref{Eq-cite1}.}

\textbf{Proof of Theorem 1}

System \eqref{Eq-1} can be rewritten in the following form:
\begin{equation}\label{Eq-2}
\begin{array}{l}
\ds{(K^2 - \varkappa^2)^{\frac{1}{4}}\, \Delta \ln (K^2 - \varkappa^2)^{\frac{1}{2}} = 4K,}\\
[2mm]
\ds{(K^2 - \varkappa^2)^{\frac{1}{4}}\, \Delta \ln \frac{K - \varkappa}{K + \varkappa} = - 4\varkappa};
\end{array}
\quad K<0,\; K^2-\varkappa^2>0.
\end{equation}

Let $K$ and $\varkappa$ be solutions of the system \eqref{Eq-2}.

We introduce the functions $\alpha$ and
$\beta$ by the formulas
\begin{equation}\label{Eq-7}
\begin{array}{l}
\ds{K= -2(\alpha^2+\beta^2),}\\
[2mm]
\ds{\varkappa} = \;\;\; 2(\alpha^2-\beta^2);
\end{array}
\quad
\alpha >0, \; \beta > 0
\end{equation}
\noindent
and system \eqref{Eq-2} takes the following form:

\begin{equation}\label{Eq-8}
\begin{array}{l}
\ds{2\sqrt{\alpha\beta}\;\Delta \ln(\alpha\beta)+8(\alpha^2+\beta^2)=0},\\
[2mm]
\ds{2\sqrt{\alpha\beta}\;\Delta \ln{\frac{\alpha}{\beta}}+4(\alpha^2-\beta^2)=0,}
\end{array}
\end{equation}
or
\begin{equation}\label{Eq-9}
\begin{array}{l}
\ds{\Delta \ln\frac{4\alpha^3}{\beta}+4\,\sqrt{\frac{4\alpha^3}{\beta}}}=0,\\
[2mm]
\ds{\Delta \ln\frac{4\beta^3}{\alpha}+4\,\sqrt{\frac{4\beta^3}{\alpha}}=0.}
\end{array}
\end{equation}

Next we put
$$p^2=\frac{4\alpha^3}{\beta}, \quad q^2=\frac{4\beta^3}{\alpha}; \quad (p>0,\; q>0).$$

Then
\begin{equation}\label{Eq-10}
K=-\frac{1}{2} \sqrt{p\,q}\,(p+q), \quad \varkappa =\frac{1}{2} \sqrt{p\,q}\,(p-q)
\end{equation}
and system \eqref{Eq-2} becomes
\begin{equation}\label{Eq-11}
\begin{array}{l}
\ds{\Delta \ln p + 2p=0,}\\
[2mm]
\ds{\Delta \ln q + 2q =0.}
\end{array}
\end{equation}

Now we apply Theorem A to any of the equations \eqref{Eq-11} and obtain

\begin{equation}\label{Eq-12}
\begin{array}{l}
\ds{p =\frac{4|w_1^\prime|^2}{(|w_1|^2+1)^2}\,,\quad w_1^\prime \neq 0,}\\
[4mm]
\ds{q =\frac{4|w_2^\prime|^2}{(|w_2|^2+1)^2}\,,\quad w_2^\prime \neq 0,}
\end{array}
\end{equation}
where $w_k$  $(k=1,2)$ are holomorphic functions in $\mathbb C$. Substituting \eqref{Eq-12} in \eqref{Eq-10}, we obtain \eqref{Eq-3}
\vskip 1mm
For the inverse, let $p$ and $q$ be two functions, given by \eqref{Eq-12}. Then equations \eqref{Eq-11} are satisfied and
the functions \eqref{Eq-3} are solutions to the system \eqref{Eq-2}. \qed
\vskip 2mm
The basic question for the system of natural PDE's \eqref{Eq-2} concerns the integrability of this system.
Theorem 1 gives explicitly the solutions of this system and therefore \eqref{Eq-2} occurs to be the first example
of an integrable system, describing a geometric class of surfaces in $\R^4$, determined by curvature conditions.

\section{When two pairs of holomorphic functions generate one and the same solution to the system of natural PDE's?}

In this section we find the relation between two pairs $(w_1, w_2)$ and $(\hat w_1, \hat w_2)$ of holomorphic functions,
which generate one and the same solution $(K, \varkappa)$ of the system \eqref{Eq-2}.

\textbf{Proof of Theorem 2}

First we shall clear up when two holomorphic functions generate one and the same solution  to the
equation \eqref{Eq-cite1}.

Let $\mathcal M: \; z=z(u,v),\; (u,v) \in \mathcal D$ be a minimal surface in $\R^3$ parameterized
by canonical principal parameters $(u,v)$ \cite{GM}.

In \cite{G} the first author proved that $\mathcal M$ has the following canonical principal representation:

\begin{equation}
z:\quad
\begin{array}{l}
\ds{z_1'=\frac{1}{2}\,\frac{w^2-1}{w'}\,,}\\
[4mm]
\ds{z_2'=-\frac{i}{2}\,\frac{w^2+1}{w'}\,,}\\
[4mm]
\ds{z_3'=-\frac{w}{w'}\,,}
\end{array}
\end{equation}
where $w$ is a holomorphic function in $\C$, which generates the solution \eqref{Eq-cite2} to \eqref{Eq-cite1}.

Putting
$$F=\frac{1}{\sqrt{-2w'}}, \quad G=\frac{w}{\sqrt{-2w'}}; \qquad \left(2FG=-\frac{w}{w'}\right)$$

we obtain the standard Weierstrass representation

\begin{equation}
z:\quad
\begin{array}{l}
z_1'=F^2-G^2\,,\\
[4mm]
z_2'=i(F^2+G^2)\,,\\
[4mm]
z_3'=2FG\,.
\end{array}
\end{equation}

Further, let $\hat{\mathcal M}$ be another minimal surface with Weierstrass representation

\begin{equation}
\hat z:\quad
\begin{array}{l}
\hat z_1'=\hat F^2-\hat G^2\,,\\
[4mm]
\hat z_2'=i(\hat F^2+\hat G^2)\,,\\
[4mm]
\hat z_3'=2\hat F \hat G\,.
\end{array}
\end{equation}

If $\mathcal M$ and $\hat{\mathcal M}$ are congruent, i.e. $\hat{\mathcal M}$ is obtained from $\mathcal M$ by a motion,
then $(F,G)$ and $(\hat F, \hat G)$ are related by the special unitary transformation (See e.g. Chapter II of \cite{N})
\begin{equation}\label{Eq-16}
\begin{array}{l}
\hat F=\;\;\; a\,F +b\,G\,,\\
[4mm]
\hat G=-\bar b\,F+\bar a \,G;
\end{array}
\quad a=\text{const}, \; b=\text{const}, \; |a|^2+|b|^2=1
\end{equation}
and vice versa.

The correspondence between the rotation and the unitary transformation is the standard spin representation of
$SO(3)$.

Replacing the expressions \eqref{Eq-16} in the formula $\ds{\hat w = \frac{\hat G}{\hat F}}$ we obtain
\begin{equation}\label{Eq-17}
\hat w = \frac{-\bar b +\bar a\,w}{a+b\,w}\,.
\end{equation}

Conversely, let $\hat w$ and $w$ be two analytic functions related by \eqref{Eq-17}.

Then
\begin{equation}
\begin{array}{l}
\ds{\hat F=\frac{1}{\sqrt{-2\hat w'}}=a\, \frac{1}{\sqrt{-2w'}}+b\,\frac{w}{\sqrt{-2w'}} = a\,F + b\, G,}\\
[4mm]
\ds{\hat G=\frac{w}{\sqrt{-2\hat w'}}=-\bar b\, \frac{1}{\sqrt{-2w'}}+\bar a\,\frac{w}{\sqrt{-2w'}} = -\bar b\,F + \bar a \, G.}
\end{array}
\end{equation}

Hence $(\hat F, \hat G)$ and $(F, G)$ are related by a special unitary transformation and
therefore $\mathcal M$ and $\hat{\mathcal M}$ are related by a motion.

Another approach to the above question has been used in \cite{OK}.

Thus we obtained the following

\begin{thm}
Let $w$ and $\hat w$ be two holomorphic functions generating one and the same
solution $\nu$ of the natural equation \eqref{Eq-cite1}.

Then
\begin{equation}
\hat w = \frac{-\bar b +\bar a\,w}{a+b\,w},\quad a={\rm const}, \; b=\rm const,\; |a|^2+|b|^2=1.
\end{equation}

Conversely, any two holomorphic functions related by \eqref{Eq-17} generate one and the same solution of \eqref{Eq-cite1}.

\end{thm}
Applying two times formula \eqref{Eq-17}, we obtain the proof of Theorem 2. \qed

\end{document}